\documentclass[twoside]{article}
\usepackage{graphics, graphicx, latexsym, amssymb, amsmath}
\newtheorem{proc}{Procedure}
\newtheorem{theo}{Theorem}
\begin{document}
\title{The perimeter generating function for \\ nondirected diagonally convex polyominoes}

\author{Svjetlan Fereti\'{c} \footnote{e-mail: svjetlan.feretic@gradri.uniri.hr} \\ 
Faculty of Civil Engineering, University of Rijeka, \\ 
Radmile Matej\v{c}i\'{c} 3, HR-51\,000 Rijeka, Croatia}

\maketitle

\begin{abstract}
A polyomino is a finite, edge-connected set of cells in the plane. At the present time, an enumeration of all polyominoes is nowhere in sight. On the other hand, there are several subsets of polyominoes for which generating functions are known. For example, there exists extensive knowledge about column-convex polyominoes, a model introduced by Temperley in 1956. While studying column-convex polyominoes, researchers also gave a look at diagonally convex polyominoes (DCPs), but noticed an awkward feature: when the last diagonal of a DCP is deleted, the remaining object is not always a polyomino. So researchers focused their attention on \textit{directed} DCPs. (A directed DCP is such a DCP that remains a polyomino when, for any $i$, its last $i$ diagonals are deleted.) Directed DCPs gradually became well understood, whereas general DCPs have remained unexplored up to now.

In this paper, we finally face \textit{general} DCPs. Modulo a little trick, which saves us from dealing with non-polyominoes, we use the layered approach (described in chapter~3 of the book ``Polygons, Polyominoes and Polycubes", edited by Anthony Guttmann). The computations are of remarkable bulk. Our main result is the perimeter generating function for DCPs; we denote it $D(d,x)$. The function $D(d,x)$ is algebraic and satisfies an equation of degree eight. The formula for $D(d,x)$ is about eight pages long. That formula involves nine polynomials in $d$ and $x$, and each of those polynomials is of degree $58$ or more in $x$. The interested reader can view the formula for $D(d,x)$ in the Maple worksheet attached to this paper.
\vspace{3mm}


\end{abstract}

\section{Introduction}

Twenty three years ago, column-convex polyominoes and directed diagonally convex polyominoes were already well explored. By contrast, there were no results about general (\textit{i.e.}, nondirected) diagonally convex polyominoes. So in her habilitation thesis, published in 1996, Mireille Bousquet-M\'{e}lou wrote: ``En comparant aux tableaux pr\'{e}c\'{e}dents les nombreux travaux effectu\'{e}s, on est d'abord frapp\'{e} par le fait que la convexit\'{e} diagonale n'a \'{e}t\'{e} que peu \'{e}tudi\'{e}e\ldots" \cite[page 52]{habilitation}. This was an open invitation to study the model. Yet, instead of new publications, what followed was just a long silence. In my case, the reason for not studying the model was the conviction that the model is unsolvable. Perhaps some other researchers thought in the same way. On the other hand, Bousquet-M\'{e}lou and Brak were more optimistic. On page 76 of the book \cite{book}, published in 2009, they wrote: ``Still, it seems that this class is sufficiently well structured to be exactly enumerable." (By ``this class" they mean nondirected diagonally convex polyominoes.) 

All in all, so far nondirected diagonally convex polyominoes have been thought of, but have not been enumerated. This state of matters now comes to an end. In the present paper, we enumerate nondirected diagonally convex polyominoes by perimeter. We start from the well-known ``turbo-Temperley" method \cite{Bousquet, Svrtan, Temperley}, also known as the layered approach, and make a little variation of it. This variation proves to suit our purposes. However, the application of (the variation of) the method requires persistent work. The computations are really cumbersome. In the end, we obtain a formula for $D(d,x)$, the perimeter generating function for nondirected diagonally convex polyominoes. (In $D(d,x)$, the exponent of $d$ is the number of diagonals and the exponent of $x$ is the perimeter.) The formula for $D(d,x)$ is eight pages long. That formula involves nine polynomials in $d$ and $x$, and each of those polynomials is of degree $58$ or more in $x$. (The highest degree is $72$.)

Let us mention that nondirected diagonally convex polyominoes have a subset with more elegant properties. That subset are \textit{directed} diagonally convex polyominoes. Let $E(d,x,y)$ be the generating function in which the coefficient of $d^k x^\ell y^m$ is the number of directed diagonally convex polyominoes with $k$ diagonals, horizontal perimeter $2\ell$ and vertical perimeter $2m$. The generating function $E$ satisfies the equation

\begin{equation}
E=d(E+1)(E+x)(E+y).
\end{equation}

\noindent Hence the number of directed diagonally convex polyominoes with $k$ diagonals is

\begin{equation}
[ d^k ]\; E(d,1,1)=\frac{1}{3k+1} {3k+1 \choose k},
\end{equation}

\noindent which is also the number of ternary trees with $k$ internal nodes. Formula (2) was found by Delest and F\'{e}dou \cite{Fedou} in 1989, and equation (1) was found by Fereti\'{c} and Svrtan \cite{Rutgers} in 1994.

This paper is organized as follows. In Section 2, we state the necessary definitions. In Section 3, we define three polyomino growing procedures. By examining those procedures, in Section 4 we establish functional equations. The functional equations are gradually solved in Section 5, and the final result is stated in Section 6. In Section 7, which is the last section, we make a comparison between our final result and the perimeter generating function for column-convex polyominoes.

\section{Definitions}

From now on, when we write diagonally convex polyomino (or DCP), we mean nondirected diagonally convex polyomino. For the standard definitions (those of a cell, of a polyomino,\ldots), we give Figure 1 as a quick reminder. More details can be found in the literature \cite{book}.

\begin{figure}
\begin{center}
\includegraphics{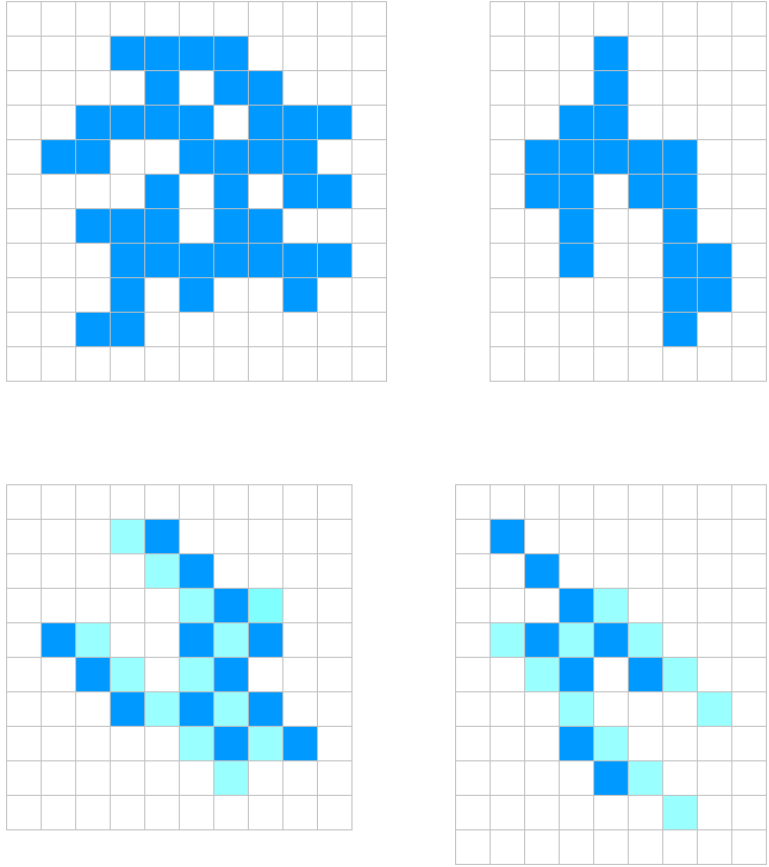} 
\caption{Top left: A polyomino. Top right: A column-convex polyomino. Bottom left: A diagonally convex polyomino. Bottom right: If diagonals were as connected as columns are, this object would be a polyomino.}
\end{center}
\end{figure}

The DCP in Figure 1 is a self-avoiding polygon, but there exist DCPs that are not self-avoiding polygons. One example is shown in Figure 2.

\begin{figure}
\begin{center}
\includegraphics{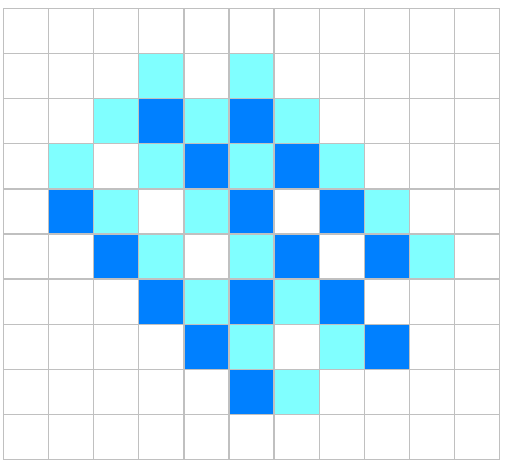} 
\caption{A diagonally convex polyomino that is not a self-avoiding polygon.}
\end{center}
\end{figure}

Given a diagonally convex polyomino $P$, we count its diagonals from southwest towards northeast. The last diagonal of $P$ is the northeasternmost diagonal of $P$.

Let $P$ be a diagonally convex polyomino with at least two diagonals. Let $a$ and $b$ be, respectively, the uppermost cell and the rightmost cell of the second-to-last diagonal of $P$. (If the second-to-last diagonal has only one cell, then $a=b$.) If the upper neighbour of $a$ is a cell of $P$, then this upper neighbour is a \textit{nose} of $P$. If the right neighbour of $b$ is a cell of $P$, then this right neighbour is a \textit{nose} of $P$. Thus, every diagonally convex polyomino has either two or one or zero noses. See Figure 3.

\begin{figure}
\begin{center}
\includegraphics{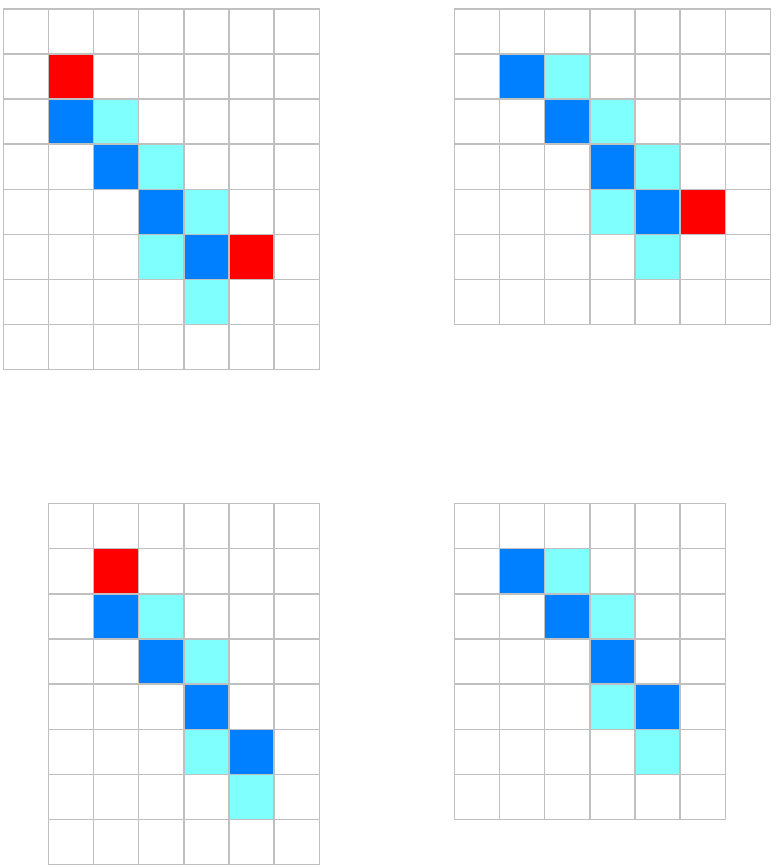} 
\caption{Four diagonally convex polyominoes. The top left one has two noses, the top right one has one nose, the bottom left one has one nose, and the bottom right one has zero noses.}
\end{center}
\end{figure}

Let $S$ be a set of diagonally convex polyominoes. We define the \textit{full generating function} of $S$ to be the formal sum

\begin{displaymath}
\sum_{P \in S} d^{di(P)} x^{pe(P)} z^{la(P)},
\end{displaymath}

\noindent where $di(P)$ stands for the number of diagonals of $P$, $pe(P)$ stands for the perimeter of $P$, and $la(P)$ stands for the size of the last diagonal of $P$.
The full generating functions for diagonally convex polyominoes with two noses, one nose and zero noses will be denoted $A(z)$, $B(z)$ and $C(z)$, respectively. (The notations $A(d,x,z)$, $B(d,x,z)$ and $C(d,x,z)$ would be somewhat lengthy.)

\section{Three polyomino growing procedures}

In the ``turbo-Temperley" method \cite{Bousquet, Svrtan, Temperley}, one first designs a procedure that adds some simple piece (typically a column) to an existing polyomino. Examining the procedure leads to functional equations, which are then solved either by iteration or by the kernel argument. We actually designed three procedures, one for each kind of diagonally convex polyominoes. Here is the first procedure:

\begin{proc}
\begin{enumerate}
\item Given a diagonally convex polyomino (say $P$) with at least two diagonals and with two noses, extend the last diagonal of $P$ by adding $r \geq 0$ new cells beyond the upper nose of $P$ and $s \geq 0$ new cells beyond the lower nose of $P$. (If $r>0$ or $s>0$, we now have an object that is not a polyomino.)
\item Add a new diagonal in such a way that the resulting object is again a polyomino and has the desired number of noses.
\end{enumerate}
\end{proc}

The first procedure is illustrated in Figure 4.

\begin{figure}
\begin{center}
\includegraphics{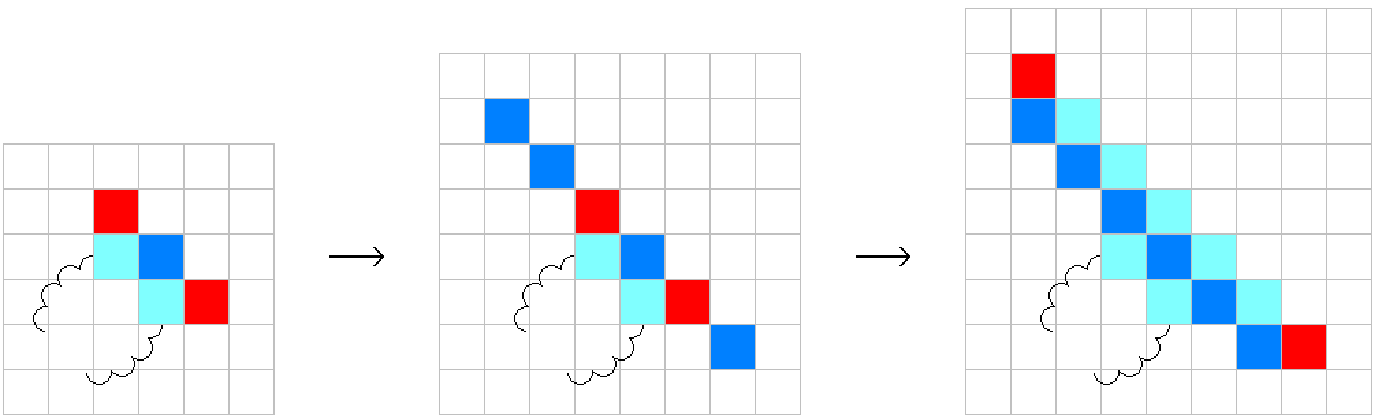} 
\caption{A polyomino with two noses produces another polyomino with two noses.}
\end{center}
\end{figure}

Here is the second procedure:

\begin{proc}
\begin{enumerate}
\item Given a diagonally convex polyomino (say $P$) with at least two diagonals and with one nose, extend the last diagonal of $P$ by adding $r \geq 0$ new cells beyond the nose of $P$. (If $r>0$, we now have an object that is not a polyomino.)
\item Add a new diagonal in such a way that the resulting object is again a polyomino and has the desired number of noses.
\end{enumerate}
\end{proc}

The second procedure is illustrated in Figure 5.

\begin{figure}
\begin{center}
\includegraphics{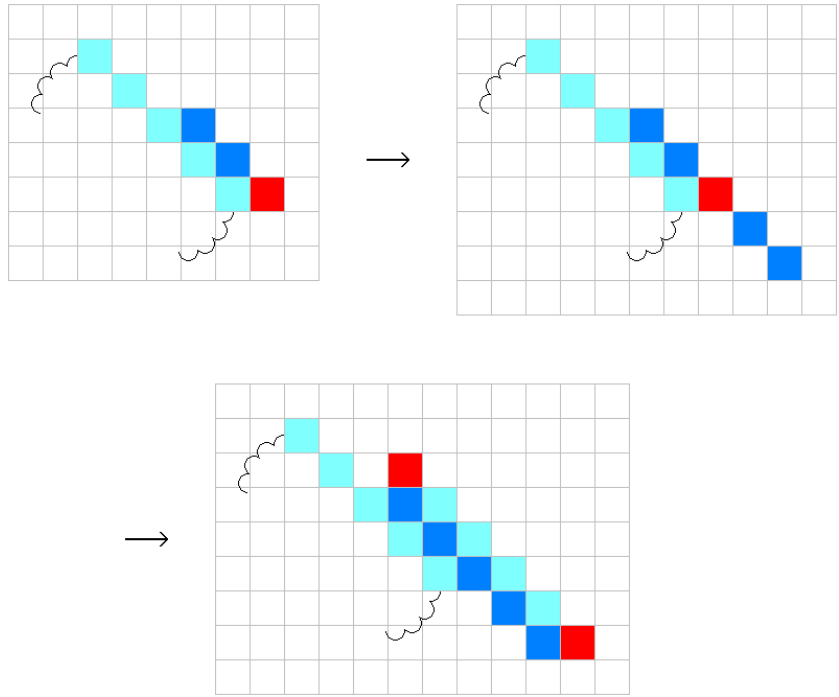} 
\caption{A polyomino with one nose produces a polyomino with two noses.}
\end{center}
\end{figure}

Here is the third (and the simplest) procedure:

\begin{proc}
Given a diagonally convex polyomino (say $P$) with at least two diagonals and with zero noses, add a new diagonal in such a way that the resulting object is a polyomino and has the desired number of noses.
\end{proc}

The third procedure is illustrated in Figure 6.

\begin{figure}
\begin{center}
\includegraphics{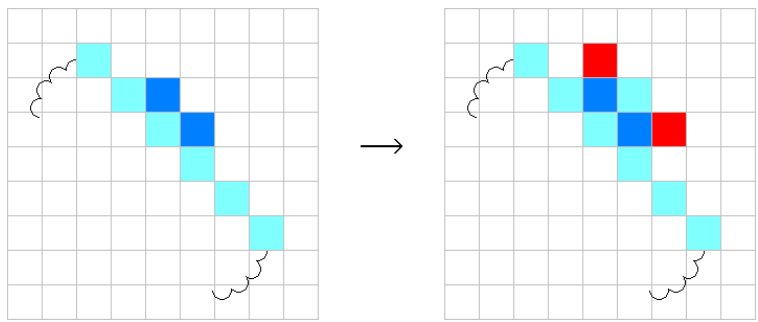} 
\caption{A polyomino with zero noses produces a polyomino with two noses.}
\end{center}
\end{figure}

Let $R$ be a diagonally convex polyomino with at least three diagonals. When we delete the last diagonal of $R$, two cases are possible. One case is that the remaining object (let us call it $Q$) is a DCP with at least two diagonals. The other case is that $Q$ has several connected components: one of those components is a DCP with at least two diagonals, and the other connected components are isolated (not-edge-contacting) cells. In either of the two cases, let $\hat {P}$ denote the connected component of $Q$ which is a DCP with at least two diagonals. Now, $\hat {P}$ is the only DCP from which the above three procedures could have produced the polyomino $R$. The polyomino $\hat {P}$ lies in the domain of only one procedure. Let us denote that procedure by $Proc$. To produce the polyomino $R$, $Proc$ first had to add the not-edge-contacting cells of $Q$ (if such cells exist). That was feasible because $\hat {P}$ has a nose wherever needed. Then $Proc$ had to add the last diagonal of $R$, which was again feasible. Thus, the above three procedures generate every DCP with at least three diagonals, and they do it in only one way.

\section{The functional equations}

Let us find a functional equation for $A(z)$. DCPs with two noses and two diagonals look like the one shown in Figure 7. If the first diagonal has $k$ cells, then the perimeter of the polyomino is $4k+4$. Thus, the part of $A(z)$ coming from polyominoes with two diagonals is 

\begin{displaymath}
\sum_{k=1}^{\infty} d^2 x^{4k+4} z^{k+1} = \frac{d^2 x^8 z^2}{1-x^4z}.
\end{displaymath}

\begin{figure}
\begin{center}
\includegraphics{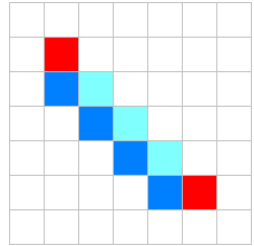} 
\caption{A DCP with two noses and two diagonals.}
\end{center}
\end{figure}

Next we consider Procedure 1. What part of $A(z)$ does come from DCPs with two noses produced from DCPs with two noses? Extending the old last diagonal, we get $4r+4s$ extra edges, and the last diagonal grows by $r+s$ cells. Then we add a new diagonal. Thus the perimeter increases by yet four edges, and the last diagonal grows by yet one cell. Therefore, this part of $A(z)$ is

\begin{displaymath}
A(z) \frac{1}{1-x^4z} \frac{1}{1-x^4z} d x^4 z = \frac{dx^4z}{(1-x^4z)^2} A(z).
\end{displaymath}

We move on to Procedure 2. What part of $A(z)$ does come from DCPs with two noses produced from DCPs with one nose? Extending the old last diagonal, we get $4r$ extra edges, and the last diagonal grows by $r$ cells. Then we add a new diagonal. Thus the perimeter increases by yet four edges, and the last diagonal grows by yet one cell. Therefore, this part of $A(z)$ is

\begin{displaymath}
B(z) \frac{1}{1-x^4z} d x^4 z = \frac{dx^4z}{1-x^4z} B(z).
\end{displaymath}

It remains to consider Procedure 3. What part of $A(z)$ does come from DCPs with two noses produced from DCPs with zero noses? When we add a new diagonal, the perimeter increases by four edges, and the last diagonal grows by one cell. Therefore, this part of $A(z)$ is

\begin{displaymath}
C(z) d x^4 z = d x^4 z C(z).
\end{displaymath}

Altogether, we have

\begin{equation}
A(z)=\frac{d^2 x^8 z^2}{1-x^4z}+\frac{dx^4z}{(1-x^4z)^2} A(z)+\frac{dx^4z}{1-x^4z} B(z)+d x^4 z C(z).
\end{equation}

In the case just considered, the desired number of noses was two. Therefore, in any of the three procedures, once the old last diagonal was extended (or not extended), the new diagonal could be added in only one way. That is why equation (3) was so easy to obtain. The functional equations for $B(z)$ and $C(z)$ cannot be obtained so quickly. Anyway, deriving the latter two functional equations is a kind of routine work (under the assumption that one knows how the functional equations are obtained in \cite{Bousquet}). Hence, we state without proof that

\begin{eqnarray}
B(z) & = & \frac{2d^2x^6z}{1-x^4z} + \frac{2dx^2z}{(1-z)(1-x^4z)}[A(1)-A(z)] + \frac{2dx^6z}{(1-x^4z)^2}A(z) \nonumber \\
& & \mbox{} + \frac{dx^2z}{(1-z)(1-x^4z)}[B(1)-B(z)] +\frac{dx^2z}{1-z}[B(1)-B(z)] +\frac{dx^6z}{1-x^4z}B(z) \nonumber \\
& & \mbox{} + \frac{2dx^2z}{1-z}[C(1)-C(z)]
\end{eqnarray}

\noindent and

\begin{eqnarray}
C(z) & = & \frac{d^2x^8z}{1-x^4z} + \frac{dz}{1-z}[A'(1)-A(1)] - \frac{dz^2}{(1-z)^2}A(1) \nonumber \\
& & \mbox{} + \frac{dz}{(1-z)^2}A(z) + \frac{2dx^4z}{(1-z)(1-x^4z)}[A(1)-A(z)] + \frac{dx^8z}{(1-x^4z)^2}A(z) \nonumber \\
& & \mbox{} + \frac{dz}{1-z}[B'(1)-B(1)] - \frac{dz^2}{(1-z)^2}B(1) + \frac{dz}{(1-z)^2}B(z) \nonumber \\
& & \mbox{} +  \frac{dx^4z}{(1-z)(1-x^4z)}[B(1)-B(z)] \nonumber \\
& & \mbox{} + \frac{dz}{1-z}[C'(1)-C(1)] - \frac{dz^2}{(1-z)^2}C(1) + \frac{dz}{(1-z)^2}C(z).
\end{eqnarray}

\section{Solution of the functional equations}

Throughout the solution, we used the computer algebra system \textit{Maple}.

\subsection{From three equations to one equation}

First we solved equations (3) and (4) for $A(z)$ and $B(z)$. Thus we got $A(z)$ and $B(z)$ expressed in terms of $A(1)$, $B(1)$, $C(1)$ and $C(z)$. Using these expressions, we eliminated $A(z)$ and $B(z)$ from  equation (5). We rearranged the equation so obtained as follows: the left-hand side was just $C(z)$; in the numerator of the right-hand side there were six unknown functions: $A(1)$, 
$A'(1)$, $B(1)$, $B'(1)$, $C(1)$ and $C'(1)$; the denominator of the right-hand side was a polynomial in $d$, $x$ and $z$.

Next we got rid of $A(1)$, $A'(1)$ and $B'(1)$. To achieve this, we set up a system of three linear equations. The first linear equation was the case $z=1$ of (3). The second linear equation was the case $z=1$ of the derivative of equation (3). (By the derivative we mean the partial derivative with respect to $z$.) The third linear equation was the case $z=1$ of equation (4). The linear system solved, we substituted its solution into the functional equation for $C(z)$. Through reducing a complex fraction to a simple fraction, the equation took the form

\begin{equation}
C(z)=\frac{P_0(d,x,z)+P_1(d,x,z)B(1)+P_2(d,x,z)C(1)+P_3(d,x,z)C'(1)}{Q(d,x,z)} \; ,
\end{equation}

\noindent where $P_0(d,x,z)$, $P_1(d,x,z)$, $P_2(d,x,z)$, $P_3(d,x,z)$ and $Q(d,x,z)$ are polynomials in $d$, $x$ and $z$.

\subsection{The roots of the denominator}

At this point, we substituted $d$ by $d^2$. (Thus we avoided dealing with non-integer powers of $d$.)

The polynomial $Q(d,x,z)$ has a nice factorization. Namely, we have

\begin{equation}
Q(d,x,z)=fa_1 \cdot fa_2 \cdot fa_3 ,
\end{equation}

\noindent where

\begin{eqnarray*}
fa_1 & = & d^2 x^{22} + x^{20} - 4d^2 x^{18} - (2d^2+5) x^{16} \\
& & \mbox{} + (d^6+2d^4+6d^2) x^{14} + (d^4+6d^2+10) x^{12} - (4d^4+4d^2) x^{10} \\
& & \mbox{} - (d^4+6d^2+10) x^8 + (2d^4+d^2) x^6 + (2d^2+5) x^4 - 1, \\
& & \\
fa_2 & = & x^4 z^2 + (d^2 x^6 - x^4 - 1) z + 1, \\
& & \\
fa_3 & = & x^8 z^4 + [-2d^2 x^{10} - (d^2+2) x^8 +2d^2 x^6 - (d^2+2) x^4] z^3 \\
& & + [d^4 x^{12} + 2d^2 x^{10} + x^8 + (4d^2+4) x^4 - 2d^2 x^2 +1] z^2 \\
& & + [-2d^2 x^6 - (d^2+2) x^4 + 2d^2 x^2 - d^2 - 2] z + 1.
\end{eqnarray*}

Factorization (7) helped us to find the roots of $Q(d,x,z)$. (By a root of $Q(d,x,z)$ we mean any $z_0$ such that $Q(d,x,z_0)=0$.) The factor $fa_1$ does not involve $z$ and hence has no roots. The equation $fa_2=0$ is quadratic; its solutions are

\begin{displaymath}
z_1=\frac{1 + x^4 - d^2 x^6 + \sqrt{1 - 2x^4 - 2d^2 x^6 + x^8 - 2d^2 x^{10} + d^4 x^{12}}}{2x^4}
\end{displaymath}

\noindent and

\begin{displaymath}
z_2=\frac{1 + x^4 - d^2 x^6 - \sqrt{1 - 2x^4 - 2d^2 x^6 + x^8 - 2d^2 x^{10} + d^4 x^{12}}}{2x^4}.
\end{displaymath}

The equation $fa_3=0$ is quartic. Quartic equations are generally not easy to handle. But fortunately, if we substitute $z$ by $x^{-2}u$, the equation $fa_3=0$ becomes symmetric: the coefficient of $u^4$ is equal to the coefficient of $u^0$, and the coefficient of $u^3$ is equal to the coefficient of $u^1$. The solution of a symmetric quartic equation boils down to the solution of two quadratic equations. Thus, we found the solutions of $fa_3=0$ rather quickly. Those solutions are

\begin{displaymath}
z_3 = \frac{H_{pl}+\sqrt{H_{pl}^2-16x^4}}{4x^4},  \qquad
z_4 = \frac{H_{pl}-\sqrt{H_{pl}^2-16x^4}}{4x^4},
\end{displaymath}

\begin{displaymath}
z_5 = \frac{H_{mi}+\sqrt{H_{mi}^2-16x^4}}{4x^4},  \qquad
z_6 = \frac{H_{mi}-\sqrt{H_{mi}^2-16x^4}}{4x^4},
\end{displaymath}

\noindent where

\begin{eqnarray*}
H_{pl} & = & 2 + d^2 - 2d^2 x^2 + 2x^4 + d^2 x^4 + 2d^2 x^6 \\
& & \mbox{} +d(1 - x^2) \cdot \sqrt{(1 + x^2)(4 + d^2 + 4x^2 - 3d^2 x^2 + 4d^2 x^4)}, \\
& & \\
H_{mi} & = & 2 + d^2 - 2d^2 x^2 + 2x^4 + d^2 x^4 + 2d^2 x^6 \\
& & \mbox{} -d(1 - x^2) \cdot \sqrt{(1 + x^2)(4 + d^2 + 4x^2 - 3d^2 x^2 + 4d^2 x^4)}.
\end{eqnarray*}

\subsection{The kernel argument}

Having found the roots of the denominator, we expanded them in power series. In the power series expansions of $z_2$, $z_4$ and $z_6$, every term has the form $c d^i x^j$, where $c \in \mathbb{Q}$, and $i,\; j \in \mathbb{N} \cup \{ 0 \}$. (This does not hold for $z_1$, $z_3$ and $z_5$. Their expansions contain some terms in which $x$ is raised to a negative power.) Now comes the kernel argument \cite{Bousquet, Svrtan}: For $i \in \{ 2,\; 4,\; 6 \}$, $C(z_i)$ is a well defined formal power series, and $Q(z_i)$ is zero. Together with (6), this means that

\begin{displaymath}
P_0(d,x,z_i)+P_1(d,x,z_i)B(1)+P_2(d,x,z_i)C(1)+P_3(d,x,z_i)C'(1)=0.
\end{displaymath}

\subsection{The remaining computations}

We now have the equations

\begin{displaymath}
P_0(d,x,z_2)+P_1(d,x,z_2)B(1)+P_2(d,x,z_2)C(1)+P_3(d,x,z_2)C'(1)=0,
\end{displaymath}

\vspace{-5mm}

\begin{displaymath}
P_0(d,x,z_4)+P_1(d,x,z_4)B(1)+P_2(d,x,z_4)C(1)+P_3(d,x,z_4)C'(1)=0,
\end{displaymath}

\noindent and
\vspace{-3mm}

\begin{displaymath}
P_0(d,x,z_6)+P_1(d,x,z_6)B(1)+P_2(d,x,z_6)C(1)+P_3(d,x,z_6)C'(1)=0.
\end{displaymath}

So, to find $B(1)$, $C(1)$ and $C'(1)$, one just has to solve a linear system. However, if one wants to find reasonably nice formulas, and not whatever formulas, then there is much more work to do. Mathematical aesthetics prevailed and we did the extra work. It consisted of several steps, of which here we shall describe the longest two. Firstly, we computed (and made use of) formulas for $z_4+z_6$, $\frac{1}{z_4}+\frac{1}{z_6}$, $z_4 z_6$ and $(z_4+z_6)^2$. The formulas for $z_4+z_6$ and $\frac{1}{z_4}+\frac{1}{z_6}$ are rather elegant:

\begin{displaymath}
z_4+z_6 = \frac{2+d^2-2d^2x^2+(2+d^2)x^4+2d^2x^6-(R_3)_{\mathrm{with}\; d^2\; \mathrm{in\; place\; of}\; d}}{2x^4},
\end{displaymath}

\begin{displaymath}
\frac{1}{z_4}+\frac{1}{z_6} = 1+\frac{d^2}{2}-d^2x^2+ \left(1+\frac{d^2}{2}\right)x^4 +d^2x^6+\frac{(R_3)_{\mathrm{with}\; d^2\; \mathrm{in\; place\; of}\; d}}{2},
\end{displaymath}

\noindent where $R_3$ is an expression defined in Theorem 1 below. The formulas for $z_4 z_6$ and $(z_4+z_6)^2$ are more complicated.

Secondly, we rationalized the lengthy denominator of $B(1)$ and $C(1)$. That denominator successfully resisted to the \textit{Maple} command \texttt{rationalize}. Hence we had to do the rationalization slowly and gradually, through writing a lot of \textit{Maple} commands.

When the formulae for $B(1)$ and $C(1)$ took a satisfactory form, the relation

\begin{equation}
A(1) = \frac{d^2 x^4 (1-x^4) [d^2 x^4 + B(1) + (1-x^4) C(1)]}{1 - (2+d^2)x^4 + x^8},
\end{equation}

\noindent gave us $A(1)$ at once. (We found relation (8) in Subsection 5.1, when we were getting rid of $A(1)$, $A'(1)$ and $B'(1)$.) At this point, diagonally convex polyominoes with two or more diagonals were all enumerated. The only diagonally convex polyomino with one diagonal is the one-celled polyomino. So we ``enumerated" the one-celled polyomino and computed the sum $d^2 x^4+A(1)+B(1)+C(1)$. Then we substituted $\sqrt{d}$ for $d$. This restored the original meaning of $d$: the exponent of $d$ was again the number of diagonals.

\section{The final result}

\begin{theo}
The perimeter generating function for diagonally convex polyominoes is given by

\begin{eqnarray}
D(d,x) & = & \frac{L+L_1 R_1+L_2 R_2+L_{1,2}R_1 R_2}{4 x^2 \Delta} \nonumber \\
& & \mbox{} -\frac{L_3 R_3+L_{1,3}R_1 R_3+L_{2,3}R_2 R_3+L_{1,2,3}R_1 R_2 R_3}{4x^2 (1-x^2) [2+d - 2dx^2 + (2+d)x^4 + 2dx^6] \Delta},
\end{eqnarray}

\noindent where

\begin{eqnarray*}
R_1 & = & \sqrt{1-2x^4-2d x^6+x^8-2d x^{10}+d^2 x^{12}}, \\
R_2 & = & \sqrt{1-(4+4d)x^2+(6+8d)x^4-(4+2d)x^6+(1-4d)x^8+2d x^{10}+d^2 x^{12}}, \\
R_3 & = & [2+4d+d^2-(4d+4d^2)x^2+(-4+6d^2)x^4+(2+4d-7d^2)x^8+(4d+4d^2)x^{10} \\
& & \mbox{} +2d^2 x^{12}+(2+4x^2+2x^4+2d x^6)R_2]^{\frac{1}{2}}.
\end{eqnarray*}

The remaining nine expressions ($\Delta$, $L$, $L_1$, $L_2$, $L_{1,2}$, $L_3$, $L_{1,3}$, $L_{2,3}$, $L_{1,2,3}$) are polynomials in $d$ and $x$. Those polynomials are too long to be stated here, but the reader can see them in the \textit{Maple} worksheet attached to this paper (and also in the Maple worksheet attached to the entry A269228 in Sloane's On-Line Encyclopedia of Integer Sequences, https://oeis.org/). The formula for $\Delta$ is shown in Figure 8 as well.
\end{theo}

\begin{figure}
\begin{center}
\includegraphics{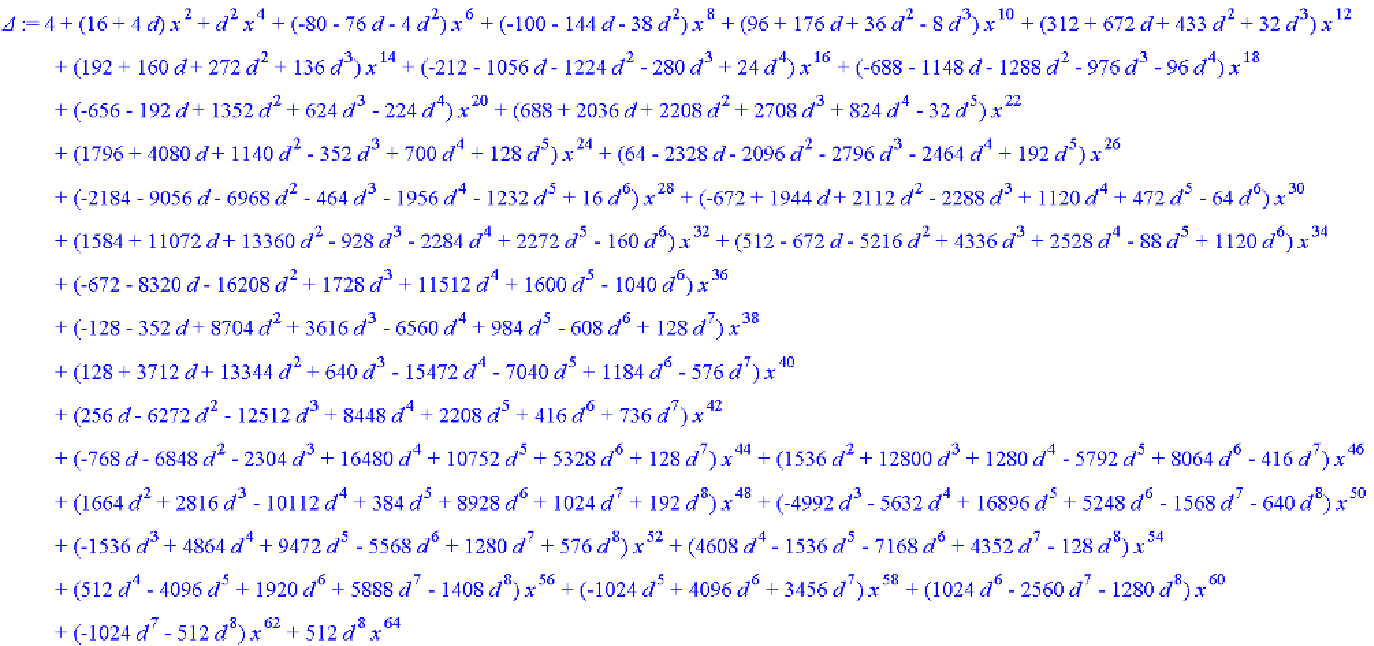} 
\caption{The formula for $\Delta$.}
\end{center}
\end{figure}

In the case $d=1$, the Taylor series expansion of $D(d,x)$ is

\begin{eqnarray}
D(1,x) & = & x^4+2x^6+7x^8+28x^{10}+122x^{12}+556x^{14}+2618x^{16}+12634x^{18} \nonumber \\
& & \mbox{}+62128x^{20}+310212x^{22}+1568495x^{24}+8014742x^{26}+41323641x^{28} \nonumber \\
& & \mbox{}+214719610x^{30}+1123244757x^{32}+5910863420x^{34}+31268459118x^{36} \nonumber \\
& & \mbox{}+166185855552x^{38}+886961294034x^{40}+\ldots
\end{eqnarray}

\section{Comparison to the perimeter generating function for column-convex polyominoes}

Row-convex polyominoes and column-convex polyominoes are essentially one and the same model. Column-convex polyominoes and diagonally convex polyominoes are two different models that show a certain degree of similarity. Here we shall take a look at the perimeter generating functions. Their structure (but not their complexity) is a matter in which diagonally convex polyominoes are not very different from column-convex polyominoes.

Let $F(x,y)$ be the power series in which the coefficient of $x^{2i}y^{2j}$ is the number of column-convex polyominoes with horizontal perimeter $2i$ and vertical perimeter $2j$. The perimeter generating function $F(x,y)$ has been studied by several authors \cite{Brak, Delest, Dubrovnik, Feretic, Svrtan, Lin, Temperley}. Modulo one wrong sign, pinpointed in \cite[page 63, footnote 1]{book}, Delest \cite{Delest} found the algebraic equation satisfied by $F(x,x)$. Then Brak \textit{et al.} \cite{Brak} computed a formula for $F(x,x)$. The next advance was Lin's \cite{Lin} formula for $F(x,y)$. That formula reads

\begin{equation}
F(x,y)=(1-y^2) \frac{\ell -6(1-y^2)r_2 -(17-x^2)(1-y^2)r_3 -r_2 r_3}{8 \delta},
\end{equation}

\noindent where

\begin{eqnarray*}
\ell & = & 42-10x^2+(-84+28x^2)y^2+(42-10x^2)y^4, \\
\delta & = & 18-2x^2+(-36+5x^2)y^2+(18-2x^2)y^4, \\
r_2 & = & \sqrt{1-2x^2+x^4+(-2-12x^2-2x^4)y^2+(1-2x^2+x^4)y^4}, \\
r_3 & = & \sqrt{2(1+x^2)(1-y^2)^2+2(1-y^2)r_2}.
\end{eqnarray*}

Later on, Svrtan and Fereti\'{c} \cite{Svrtan} found two more formulas for $F(x,y)$. Those two formulas are

\begin{equation}
F(x,y)=(1-y^2)\left[1-\frac{2\sqrt{2}}{3\sqrt{2}-\sqrt{1+x^2+\sqrt{(1-x^2)^2-\frac{16x^2y^2}{(1-y^2)^2}}}}\right]
\end{equation}

\noindent and

\begin{equation}
F(x,y)=(1-y^2)\left[1-\frac{4}{6-\sqrt{(1-x)^2-\frac{4xy^2}{1-y^2}}-\sqrt{(1+x)^2+\frac{4xy^2}{1-y^2}}}\right].
\end{equation}

In the long shadow of the formula for $D(d,x)$, formulas (11), (12) and (13) all look simple. However, taking a closer look, we see that formulas (12) and (13) are shorter than (11). We wrote (12) and (13) in one line, but for (11) we needed several lines. It is plausible that the formula for $D(d,x)$ also can be written in a reasonably shorter way. Perhaps that could be achieved via some multiparametric enumeration of column-convex polyominoes on the hexagonal lattice. Namely, a diagonally convex polyomino can be viewed as a polyomino with hexagonal cells which, in addition to being column-convex, has the property that every pair of its cells is connected by a path that does not make vertical steps. See Figure 9.

\begin{figure}
\begin{center}
\includegraphics{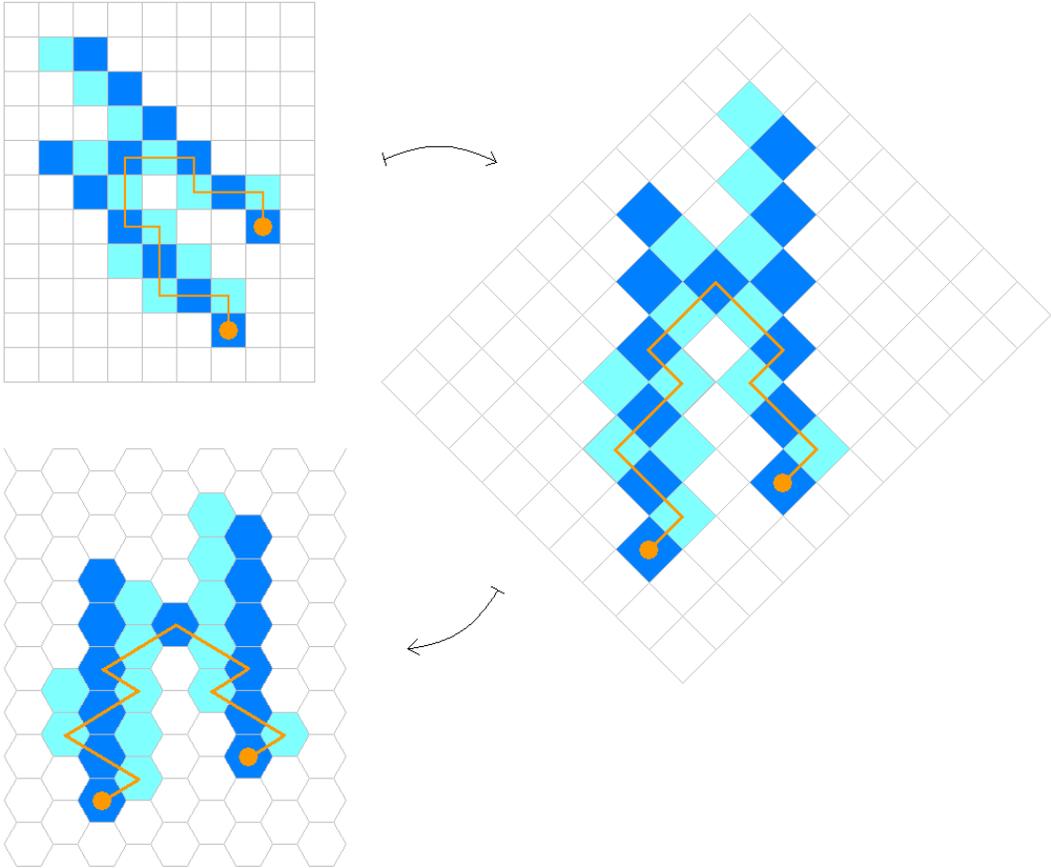} 
\caption{The transformation of a diagonally convex polyomino into a column-convex polyomino with hexagonal cells.}
\end{center}
\end{figure}

A good point of expression (11) is that it is in a standard form. So we can readily compare our formula (9) to (11). The most interesting part of formulas (9) and (11) are the expressions that involve square roots. In (9), there are three such expressions: $R_1$, $R_2$ and $R_3$. In (11), there are two such expressions: $r_2$ and $r_3$. True to their names, $R_2$ and $R_3$ are counterparts of $r_2$ and $r_3$, respectively. On the other hand, the expression $R_1$ is endemic. In the perimeter generating function of column-convex polyominoes, there is no expression of that kind. Hence the perimeter generating function $D(d,x)$ satisfies an octic equation, while the perimeter generating function $F(x,y)$ satisfies a quartic equation.

The Taylor series expansion of $F(x,x)$ is

\begin{eqnarray}
F(x,x) & = & x^4+2x^6+7x^8+28x^{10}+122x^{12}+558x^{14}+2641x^{16}+12822x^{18} \nonumber \\
& & \mbox{}+63501x^{20}+319554x^{22}+1629321x^{24}+8399092x^{26}+43701735x^{28} \nonumber \\
& & \mbox{}+229211236x^{30}+1210561517x^{32}+6432491192x^{34}+34364148528x^{36} \nonumber \\
& & \mbox{}+184463064936x^{38}+994430028087x^{40}+\ldots
\end{eqnarray}

Let $cc_n$ denote the number of column-convex polyominoes with perimeter $n$, let $dc_n$ denote the number of diagonally convex polyominoes with perimeter $n$, and let $ra_n=\frac{cc_n}{dc_n}$. From the Taylor series expansions (10) and (14) we see that $ra_{4}=1$, $ra_{6}=1$,..., $ra_{12}=1$. Then the ratios start to grow, but the growth rate is low: we have $ra_{14}=1.0036$, $ra_{16}=1.0088$,..., $ra_{36}=1.0990$, $ra_{38}=1.1100$, $ra_{40}=1.1212$. To get more evidence, we computed the Taylor series expansions up to the terms with $x^{200}$. It turned out that the trend continues. For example, we have $ra_{196}=2.5646$, $ra_{198}=2.5922$, and $ra_{200}=2.6201$. Thus, if diagonally convex polyominoes were a subset of self-avoiding polygons, they would approximate all self-avoiding polygons almost as closely as column-convex polyominoes do. (This claim is based on numerical evidence; we have not done the asymptotic analysis.) However, diagonally convex polyominoes \textit{are not} a subset of self-avoiding polygons. Recall Figure 2.

\end{document}